\documentclass[12pt]{amsart}
\usepackage{amssymb,amsmath,amsthm}

\newtheorem{thm}{Theorem}

\begin{document}
\pagestyle{myheadings}
\title
{Orthogonal polynomials of discrete variable and boundedness of
Dirichlet kernel}
\author{Josef Obermaier and Ryszard Szwarc}
\address{J. Obermaier, Institute of Biomathematics and Biometry,
GSF-National Research Center for Environment and Health,
Ingolst\"adter Landstrasse 1,
D-85764 Neuherberg, Germany}
\email{josef.obermaier@gsf.de}
\address{R. Szwarc, Institute of Mathematics,
University of Wroc\l aw, pl.\ Grunwaldzki 2/4, 50-384 Wroc\l aw, Poland}
\email{szwarc@math.uni.wroc.pl}
\date{}
\thanks{Supported by European Commission Marie Curie Host
Fellowship for the Transfer of Knowledge ``Harmonic Analysis, Nonlinear
Analysis and Probability''  MTKD-CT-2004-013389,  KBN (Poland) under Grant 2 P03A 028 25 and
DFG Contract 436 POL 17/1/04.}
\begin{abstract}
For orthogonal polynomials defined by compact Jacobi matrix with
exponential decay of the coefficients, precise properties of
orthogonality measure is determined. This allows showing uniform
boundedness of partial sums of orthogonal expansions   with
respect to $L^\infty$ norm, which generalize analogous results
obtained for little $q$-Legendre, little $q$-Jacobi and little
$q$-Laguerre polynomials, by the authors.
\end{abstract}
\maketitle
\markboth{\normalsize\sl  Ryszard Szwarc and Josef Obermaier}{{\normalsize\sl  Orthogonal polynomials of discrete variable}}
\section{Introduction}
Let $s_n(f)$ denote the $n$th partial sum of the classical Fourier
series of a continuous $2\pi$ periodic function $f(\theta).$ We
know that the quantities $\|s_n(f)\|_\infty$ need not to be
uniformly bounded since the Lebesgue numbers
$\int_0^{2\pi}|D_n(\theta)|d\theta$ behave like constant multiple
of $\log n,$ where $D_n$ denotes the Dirichlet kernel.

In principle this is Faber's result \cite{F} which shows that the
system of trigonometric polynomials does not constitute a Schauder
basis with respect to the set of continuous functions
$C([0,2\pi])$. Moreover in case of $C([-1,1])$ he derived the
analogous result regarding a system of algebraic polynomials with
degrees increasingly passing through all positive integers. Let us recall
that a sequence $\{\varphi_n\}_{n=0}^{\infty}$ in $C(S)$, where
$S\subset \mathbb R$, is called a Schauder basis with respect to
$C(S)$ if for every $f\in C(S)$ there exists a unique sequence of
numbers $\{a_n\}_{n=0}^{\infty}$ such that
\begin{equation}
f=\sum_{n=0}^{\infty}a_n \varphi_n.
\end{equation}
Privalov \cite{P1} refined the result of Faber: If
$\{P_n\}_{n=0}^\infty$ is a Schauder basis with respect to
$C([a,b])$ consisting of algebraic polynomials then there are
$\epsilon >0$ and $m \in {\mathbb N}_0$ such that $\deg P_n\ge
(1+\epsilon) n$ for all $n\ge m$. Other way around, he proved in
\cite{P2} a remarkable result that for any $\epsilon >0$ there exists an algebraic
polynomial Schauder Basis $\{P_n\}_{n=0}^\infty$ with $\deg P_n
\le (1+\epsilon) n.$ Such a basis is called
basis of optimal degree with respect to $\epsilon$. Concerning
the existence of an orthogonal polynomial Schauder basis of
optimal degree there are two particular results we want to mention.
The first gives orthogonal basis with respect to Tchebyshev weight of the first kind \cite{KPS}
and the second   with respect to the Legendre weight \cite{SK}.
The problem of construction or even of the existence of a
minimal basis for general Jacobi weights seems still to be open
and, more generally, it is open for an arbitrary positive measure
concentrated on an interval.

There are reasons for to have a polynomial basis
$\{P_n\}_{n=0}^\infty$ with $\deg P_n = n$. For instance this
would imply that the partial sums $s_n(f)$ are converging towards
$f$ with the same order of magnitude as the elements of best
approximation in ${\mathcal P}_n$ do \cite[§ 19, Theorem
19.1]{SI}, where ${\mathcal P}_n$ denotes the set of algebraic
polynomials with degree less than or equal to $n$. With this in mind and
due to the results above, we have to switch to spaces $C(S)$,
where $S$ differs from an interval.

The question arises: Do there exist a measure space and a
corresponding orthogonal polynomial system $\{R_n\}_{n=0}^\infty$
with $\deg R_n = n$ such that the partial sums of the Fourier
series are uniformly bounded in $\|\cdot\|_\infty$ norm ?

The situation is trivial if the support is finite. But the
problem becomes nontrivial if the measure space is infinite, for
instance of the form $\{q^n\}_{n=0}^\infty$ for some number
$0<q<1.$ There are examples of systems of orthogonal polynomials
whose orthogonality
 measure is concentrated on the sequence $\{q^n\}_{n=0}^\infty.$ Little $q$-Legendre polynomials,
 more generally $q$-Jacobi polynomials and
 little $q$-Laguerre polynomials are such. The uniform boundedness of $\|s_n(f)\|_\infty$ have been shown
 for these systems in \cite{O, OS1, OS2}. The proof depended heavily on the precise knowledge
 of the orthogonality measure and pointwise estimates of these polynomials.

 In this paper we will generalize by far these results by allowing general orthogonal polynomials
 satisfying a three term recurrence relation
 $$xp_n=-\lambda_np_{n+1}+\beta_np_n-\lambda_{n-1}p_{n-1},$$
 where $\lambda_n>0,\ \beta_n\in \mathbb{R}.$ Assuming boundedness
 of these coefficients the orthogonality measure $\mu$ on the real line is determined uniquely.
 However finding this measure explicitly is a hopeless task in general and can be achieved
 in very few special cases. Nonetheless we are able sometimes to derive certain properties
 of this measure. We will use the well known fact that if $J$ is the Jacobi matrix associated with
 the coefficients $\{\lambda_n\}_{n=0}^\infty$ and $\{\beta_n\}_{n=0}^\infty,$ i.e.
\begin{equation}\label{jacobi}
J=\begin{pmatrix}
\beta_0 & \lambda_1 & 0 & 0 & \cdots \\
\lambda_1 & \beta_1 & \lambda_2 & 0 & \cdots \vspace{-3pt}\\
0 & \lambda_2 & \beta_2 & \lambda_3 & \ddots \vspace{-3pt}\\
0 & 0 & \lambda_3 & \beta_3 & \ddots \\
\vdots & \vdots & \ddots & \ddots & \ddots
\end{pmatrix},
\end{equation}
then the spectrum of $J$ on $\ell^2( \mathbb{N}_0)$ coincides with the support of $\mu.$

 In this paper we impose conditions on the sequences $\{\lambda_n\}_{n=0}^\infty$ and
 $\{\beta_n\}_{n=0}^\infty$ so that determining the behavior of the orthogonality measure
 is possible. In particular we will assume that these coefficients
 have exponential decay at infinity. The properties of the orthogonality measure
 will be sufficient for proving
 the uniform boundedness of the norms $\|s_n\|_{L^\infty\to L^\infty}.$

Throughout the paper we will be using certain classical results
concerning orthogonal polynomials. In most such cases   references
will be given. In particular we will use the following well known
property, whose proof follows immediately from orthogonality. If
$\mu((a,b))=0,$ where $\mu$ is an orthogonality measure, then the
polynomial $p_n$ may have at most one root in the interval
$[a,b].$ Moreover, if $\mu((c,+\infty))=\mu((-\infty,d))=0$ then
$p_n$ does not vanish in either interval.

By $a_n\approx b_n$ we will mean that the ratio $a_n/b_n$ has a positive limit, while by
$a_n\sim b_n$ we will mean that the ratio $a_n/b_n$ is positive, bounded and bounded away from zero.

{\bf Acknowledgment} We thank Walter Van Assche for turning our
attention to Tchebyshev-Markov-Stieltjes inequalities.
\section{Orthogonality measure}
Let $R_n(x)$ denote polynomials satisfying  a three term recurrence relation
\begin{equation}\label{basic}
xR_n(x)=-\gamma_nR_{n+1}(x)+\beta_nR_n(x)-\alpha_nR_{n-1}(x),
\end{equation}
where $\alpha_0=0$ and $R_0(x)\equiv 1.$ We assume that $\gamma_n, \alpha_{n+1}>0$ and
\begin{equation}\label{beta}
\beta_n=\alpha_n+\gamma_n.\end{equation} In this way the polynomials are normalized at $0$ so that
\begin{equation}\label{norm}
R_n(0)=1.
\end{equation}
Since the coefficient of the leading term  of $R_n$ is alternating,
and  the roots of $R_n$ are distinct and real (see \cite[Theorem I.5.2]{AH}),
all these roots are positive in view of (\ref{norm}).
Therefore (see \cite[Proof of Theorem 2.1.1, for $\tau=0$]{AH}) there is an orthogonality measure $\mu$ supported on half line $[0,+\infty).$
Let $h(0)=1$ and
$$ h(n)={\gamma_0\gamma_1\ldots \gamma_{n-1}\over \alpha_1\alpha_2\ldots \alpha_n}.
$$
It can be easily computed that the polynomials
 \begin{equation}\label{orth} p_n(x)=\sqrt{h(n)}R_n(x)
 \end{equation}
are orthonormal and satisfy the recurrence relation
\begin{equation}
xp_n(x)=-\lambda_np_{n+1}(x)+\beta_np_n(x)-\lambda_{n-1}p_{n-1}(x),\label{recur}
\end{equation}
where
\begin{equation}\label{lambda}
\lambda_n=\sqrt{\alpha_{n+1}\gamma_n}.
\end{equation}

We will consider polynomials with special properties such that
the orthogonality measure is concentrated on a sequence of points $\xi_n$ such that
$\xi_n \searrow 0$ when $n\to \infty.$  There are many instances of such behavior, e.g.
little $q$-Jacobi polynomials, little $q$-Laguerre polynomials. Also we require that the polynomials satisfy
nonnegative product linearization  property, i.e.
the coefficients in the expansions
\begin{equation}
R_n(x)R_m(x)=\sum_{k=|n-m|}^{n+m}g(n,m,k)R_k(x)
\end{equation}
are all nonnegative. The above mentioned polynomials fulfill this property for certain values of
parameters.

We will deal with general orthogonal polynomials satisfying the two above properties.
In order to ensure
the proper behaviour of the orthogonality measure as well as nonnegative linearization property
we assume that there are   constants $q,$  $\kappa,$ $s,$
 $c$ and $N$ such that

\begin{align}\label{compact}
&\alpha_n\approx   q^n,\ \gamma_n\approx   q^n,&&0<q<1, \\
&\alpha_n \le \kappa\gamma_n, && 1\le \kappa<q^{-1}+q-1, \label{tech}\\
&h(n)\sim s^n  && s>1,\label{s}\\
& \lambda_n\le \beta_{n+1}-c\beta_{n+2},&& {1+q\over 1+q^2}<c<{1\over q}.\label{ad}\\
& \beta_1\le \beta_0.&&\label{last}\\
& \beta_n-c\beta_{n+1}\ge \beta_{n+1}-c\beta_{n+2},\quad n\ge N. &&\label{last1}
\end{align}

{\bf Remark.} Assumption (\ref{tech}) is technical. In many cases, like little $q$-Jacobi polynomials,
this assumption is satisfied with $\kappa=1.$ Actually it is natural to expect $\alpha_n\le \gamma_n$
(see (\ref{hn})).

By assumptions (\ref{ad}) and (\ref{last}) we obtain that $\beta_n$ is a decreasing sequence and
\begin{equation}\label{msz}
\lambda_n\le \beta_{n+1}-\beta_{n+2}, \quad n\ge 0.
\end{equation}
Hence the assumptions  of  \cite[Theorem 1]{MS} are satisfied. The fact that $\beta_n$ is decreasing
instead of being increasing   follows from normalizing our polynomials
 in such a way that the sign of the leading coefficient is alternating, instead of
 being positive like in \cite{MS}.
Therefore the polynomials
$\{R_n\}_{n=0}^\infty$ admit nonnegative product linearization.
This property implies that (see \cite[(17), p. 166]{sch})
$$
|R_n(x)|\le 1, \qquad x\in {\rm supp \mu},
$$
or equivalently
\begin{equation}\label{B}
|p_n(x)|\le p_n(0), \qquad x\in {\rm supp \mu}.
\end{equation}
By orthonormality and by (\ref{B}) we have $p_n^2(0)\ge 1.$
In particular
\begin{equation}\label{hn}
h(n)=p_n^2(0)={\gamma_0\gamma_1\ldots \gamma_{n-1}\over \alpha_1\alpha_2\ldots \alpha_n}\ge 1.
\end{equation}
In the next theorem we are going to describe the orthogonality measure $\mu$ for the polynomials $\{R_n\}_{n=0}^\infty.$
\begin{thm}
Assume (\ref{basic}-\ref{last1}) are satisfied.
Then the orthogonality measure $\mu$ is concentrated on decreasing sequence
$\{\xi_n\}_{n=1}^\infty,$ where
$\xi_n\sim q^n,$ the quantity $1-\xi_{n+1}/\xi_n$ is bounded away from 0, and $\mu([0,\xi_n])\sim s^{-n}.$
\end{thm}
\noindent{\bf Remark 1.} The conclusion of the theorem cannot be
strengthened to $\mu(\{\xi_n\})\sim s^n.$ Indeed, consider the
probability measure
$$\mu={3\over 2}\sum_{n=0}^\infty {1\over 4^{n+1}}\delta_{2^{-2n}}+ {7\over 2}\sum_{n=0}^\infty {1\over 8^{n+1}}\delta_{2^{-(2n+1)}}.$$
Then $\xi_n\sim 2^{-n}$ and $\mu([0,\xi_n])\sim 2^{-n}$ but
$\mu(\{\xi_n\})\not\sim 2^{-n}.$ Of course we cannot guarantee  that the polynomials orthogonal with respect
to this measure satisfy nonnegative product linearization.
\begin{proof}
Let $J$ denote the Jacobi matrix associated with the polynomials $p_n$ (see (\ref{jacobi})).
 By assumptions $J$ is a   compact operator on $\ell^2( \mathbb{N}_0).$ Moreover
 $J$ is semipositive definite because by (\ref{lambda}) we have
$J=S^*S,$ where
$$S
=\begin{pmatrix}
\sqrt{\gamma_0} & \sqrt{\alpha_1} & 0 & 0 & \cdots \\
0& \sqrt{\gamma_1} & \sqrt{\alpha_2} & 0 & \cdots \vspace{-3pt}\\
0 & 0 & \sqrt{\gamma_2} & \sqrt{\alpha_3} & \ddots \vspace{-3pt}\\
0 & 0 & 0 & \sqrt{\gamma_3} & \ddots \\
\vdots & \vdots & \ddots & \ddots & \ddots
\end{pmatrix}.
$$
Hence the  spectrum of $J$
consists of $0$ and a decreasing sequence of points $\{\xi_n\}_{n=1}^\infty$ accumulating at zero.
As we have mentioned in the introduction, the support of $\mu$ coincides with the spectrum of $J.$ First we will show that
$\mu(\{0\})=0.$   Indeed, by \cite[Theorem 2.5.3]{AH} we have
$$\mu(\{0\})^{-1}=\sum_{n=0}^\infty p_n^2(0).$$
We know that $p_n^2(0)\ge 1$ (see (\ref{hn})). Hence $\mu(\{0\})^{-1}=\infty.$

Now we turn to determining the behavior of $\xi_n.$ Let $\{x_{jn}\}_{j=1}^n$ denote the zeros of the  polynomial
$p_n(x)$ arranged in the increasing order. It is well known (see \cite[Exercise I.4.12]{CH}) that this set coincides with the set
of eigenvalues of the truncated Jacobi matrix $J_n,$ where
$$
J_n=
\begin{pmatrix}
\beta_0 & {\lambda_0}& 0 &\cdots &0 &0 \\
{\lambda_0} &\beta_1 & {\lambda_1}&\cdots &0& 0\\
0&{\lambda_1} & \beta_{2} &\cdots &0& 0\\
\vdots & \vdots & \vdots &  \ddots &\vdots & \vdots\\
0   & 0 & 0 &\cdots & \beta_{n-2} & {\lambda_{n-2}}\\
0   & 0 & 0 &\cdots &  {\lambda_{n-2}}& \beta_{n-1}\\
\end{pmatrix}.
$$
By (\ref{ad}) and by (\ref{last}) we have for $n\ge 2$
\begin{eqnarray*}
&&\lambda_0\le \beta_1-c\beta_2\le \beta_0-c\beta_n,\\
&&\lambda_{i-1}+\lambda_i\le \beta_i-c\beta_{i+2}\le \beta_i-c\beta_n,\quad 1\le i\le n-2,\\
&&\lambda_{n-2}\le \beta_{n-1}-c\beta_n.
\end{eqnarray*}
These inequalities imply
$$J_n\ge c\beta_nI_n,$$
where
$I_n$ denotes the identity matrix of rank $n.$
  Therefore $x_{1n}\ge c\beta_n.$ On the other hand by orthogonality
the polynomial $p_n(x)$ cannot change
sign more than once between two consecutive points of ${\rm supp\ }\mu$ and it cannot change sign in the interval
$[\xi_1,+\infty).$ Therefore $\xi_n\ge x_{1n}$ and consequently
\begin{equation}\label{lower}
\xi_n\ge c\beta_n.
\end{equation}

For the upper estimate we will use the minimax theorem. Let $(\cdot,\cdot)$ denote the standard inner product
in the real Hilbert space $\ell^2( \mathbb{N}_0)$ and  $\{\delta_n\}_{n=0}^\infty$  denote the standard orthogonal basis in this space.
We have
$$
\xi_{n}=\min_{v_1,\ldots,v_{n-1}}\max_{w\perp v_1,\ldots,v_{n-1}}{(Jw,w)\over (w,w)}\le
\max_{w\perp \delta_0,\ldots, \delta_{n-2}}{(Jw,w)\over (w,w)}=\|A_n\|,
$$
where
$$
A_n=\begin{pmatrix}
\beta_{n-1} & \lambda_{n-1} & 0 & 0 & \cdots \\
\lambda_{n-1} & \beta_n & \lambda_n & 0 & \cdots \vspace{-3pt}\\
0 & \lambda_n & \beta_{n+1} & \lambda_{n+1} & \ddots \vspace{-3pt}\\
0 & 0 & \lambda_{n+1} & \beta_{n+2} & \ddots \\
\vdots & \vdots & \ddots & \ddots & \ddots
\end{pmatrix}.
$$
Therefore
$$\|A_n\|\le \max\{\beta_{n-1}+\lambda_{n-1}, \max\{\lambda_{i-1}+\beta_i +\lambda_i\,:\, i\ge n\}\}.$$
By   (\ref{ad})   we obtain
\begin{eqnarray*}
 \beta_{n-1}+\lambda_{n-1}&\le &\beta_{n-1}+\beta_n-c\beta_{n+1},\\
 \lambda_{i-1}+\beta_i +\lambda_i&\le &2\beta_i-c\beta_{i+1},\qquad\qquad i\ge n.
\end{eqnarray*}
By (\ref{last1}) and the fact that $\beta_n$ is decreasing we may conclude that
$$ \|A_n\|\le \beta_{n-1}+\beta_n-c\beta_{n+1},$$
for $n\ge N.$

Summarizing we proved that
\begin{equation}\label{xi}
c\beta_n\le \xi_n\le \beta_{n-1}+\beta_n-c\beta_{n+1}, \quad n\ge N,
\end{equation}
which shows that $\xi_n\sim q^n$ because $\beta_n=\alpha_n+\gamma_n\approx q^n.$
For $n\ge N$ we have
$$ \xi_{n+1}\le \beta_{n}+\beta_{n+1}-c\beta_{n+2}.$$
Thus
$$
{\xi_{n+1}\over \xi_{n}}\le {\beta_{n}+\beta_{n+1}-c\beta_{n+2}\over c\beta_n}.
$$
By $\beta_n\approx q^n$ and by the second part of (\ref{ad}) we
obtain
$$\limsup_{n\to \infty}{\xi_{n+1}\over \xi_{n}}={1+q-cq^2\over c}<1.$$
Since $\xi_{n+1}<\xi_n$ for any $n,$ the quantity $\xi_{n+1}/\xi_n$ is bounded away
from zero.

Concerning the second part we will estimate from above the quantities
$$ \mu(\{\xi_n\})=\left (\sum_{j=0}^\infty p_j(\xi_n)^2\right )^{-1}.$$
There is a  positive constant $C $ such that
\begin{equation}
 {\xi_n\over \gamma_j}\le Cq^{n-j}.\label{qn}
\end{equation} By the second part of (\ref{tech}) there exists a positive integer $t$  such that
\begin{equation}
q^t\le {1-q\over C}\left (1-{\kappa q\over 1-q+q^2}\right )\label{cond}
 \end{equation}
We are going to show   that for $j\le n-t$ there holds $R_{j-1}(\xi_n)> 0$ and
\begin{eqnarray}
{R_{j}(\xi_n)\over
R_{j-1}(\xi_n)}=1-\varepsilon_j,\label{quo}\\
0\le \varepsilon_j \le {Cq^{n-j}\over 1-q}.\label{q}
\end{eqnarray}
The proof will go by induction on $j\le n-t.$
By (\ref{basic}) and by (\ref{qn}) we have for $j=1$
$${R_{1}(\xi_n)\over
R_{0}(\xi_n)}=R_{1}(\xi_n)=1-{\xi_n\over \gamma_0},$$
and $$\varepsilon_1={\xi_n\over \gamma_0}\le Cq^n\le {Cq^{n-1}\over 1-q}.$$
Assume that (\ref{quo}) and  (\ref{q}) hold for  $j,$ where $0 \le j<n-t.$
Hence by (\ref{cond}) we obtain
$$\varepsilon_j\le {Cq^{n-j}\over 1-q}\le {Cq^t\over 1-q}\le 1-
{\kappa q\over q^2-q+1}<1,
$$ which by (\ref{quo})
implies $R_{j}(\xi_n)> 0.$
By virtue of (\ref{basic}) and $\beta_j=\alpha_j+\gamma_j$ we have
$$\gamma_j{R_{j+1}(\xi_n)\over R_{j}(\xi_n)}+\alpha_j{R_{j-1}(\xi_n)\over R_{j}(\xi_n)}=
\alpha_j+\gamma_j-\xi_n.$$
Therefore
\begin{equation}\label{imp}
\varepsilon_{j+1}={\xi_n\over \gamma_j}+{\alpha_j\over \gamma_j}{\varepsilon_j\over 1-\varepsilon_j}.
\end{equation}
By induction hypothesis, in view of  (\ref{tech}) and (\ref{qn}),   we get
\begin{multline*}
\varepsilon_{j+1}\le Cq^{n-j}+\kappa {Cq^{n-j}\over
1-q-Cq^{n-j}} \\=
{Cq^{n-j-1}\over 1-q} \,\left [q(1-q)+ {\kappa q(1-q)\over 1-q - Cq^{n-j}}\right ]\\
\le {Cq^{n-j-1}\over 1-q}\,\left [q(1-q)+ {\kappa q(1-q)\over 1-q - Cq^{t}}\right ],
\end{multline*}
because $n-j\ge t.$
Condition (\ref{cond})
%and the fact that $\kappa < q^{-1}+q-1$
implies
$$1-q-Cq^t \ge {\kappa q(1-q)\over 1-q+q^2}.$$
Therefore
$$ q(1-q)+ {\kappa q(1-q)\over 1-q - Cq^{t}}\le q(1-q)+1-q+q^2 = 1.$$
Therefore
$$\varepsilon_{j+1}
\le {Cq^{n-j-1}\over 1-q}.
$$
The assumption (\ref{cond}) and $\kappa\ge 1$ imply
$$q^t\le {1\over C}(1-q)^2.$$
Now   (\ref{quo}) and (\ref{q})   yield that for $j\le n-t$ there holds
\begin{multline*}
R_j(\xi_n)=(1-\varepsilon_1)(1-\varepsilon_1)\ldots (1-\varepsilon_j)\ge
1-\sum_{i=1}^j \varepsilon_i\\ \ge
1-\sum_{i=1}^j {Cq^{n-i}\over 1-q}\ge 1 - {Cq^{n-j}\over (1-q)^2}
\ge 1-{Cq^t\over (1-q)^2}>0.
\end{multline*}
Let $\eta=1-{Cq^t\over (1-q)^2}.$ Then $R_j(\xi_n)\ge \eta,$ for $j\le n-t.$
In view of $R_j(x)=p_j(x)/p_j(0)$ we get
$$
p_j(\xi_n)\ge \eta p_j(0), \qquad j\le n-t.
$$
Therefore
$$\mu(\xi_n)^{-1}=\sum_{j=0}^\infty p_j^2(\xi_n)
\ge \sum_{j=0}^{n-t} p_j^2(\xi_n) \ge \eta^2\sum_{j=0}^{n-t}p_j(0)^2.$$
By (\ref{s}) and (\ref{hn}) we have $p_j^2(0)\sim s^j$ for $s>1.$ Hence
$$\mu(\xi_n)\le Ds^{-n}$$ for some constant $D.$
This implies
$$\mu([0,\xi_n])=\mu((0,\xi_n])=\sum_{k=n}^\infty \mu(\xi_k)\le {D\over s-1}\,s^{-n-1}.$$

It remains to show that $\mu([0,\xi_n])\ge ds^{-n}$ for some constant $d.$
To this end we will use Tchebyshev inequalities. Let $\{x_{ni}\}_{i=1}^n$ denote the zeros
of the polynomial $p_n$ arranged in the increasing order. Let
$$\mu_{ni}= \left(\sum_{j=0}^{n-1}p_j^2(x_{ni})\right)^{-1}.$$
By \cite[Thm. 3.41.1]{sz} we have
$$ \mu_{n1}\le \mu ([0,x_{n2})).$$
Since $|p_j(x_{n1})|\le p_j(0)$ (see (\ref{B})) we have
$$\mu ([0,x_{n2}))\ge \left(\sum_{j=0}^{n-1}p_j^2(0)\right)^{-1}\ge ds^{-n}$$
for some $d>0.$ By orthogonality no two consecutive points of $\{x_{ni}\}_{i=1}^n$ may lie between two consecutive
points of $\{\xi_m\}_{m=1}^\infty.$ Also $x_{nn}< \xi_1.$ Therefore $x_{n2}< \xi_{n-1}.$ This gives
$$\mu([0,\xi_n])=\mu ([0,\xi_{n-1}))\ge \mu ([0,x_{n2}))\ge ds^{-n}.$$\end{proof}
\section{Boundedness of Dirichlet kernel}
Consider orthogonal polynomials defined by (\ref{basic}).
 Let $\mu$ denote the corresponding orthogonality measure.
Let $S=\mbox{supp}\,\mu.$

For functions $f\in C(S)$ and $k\in \mathbb{N}_0$ the generalized
Fourier coefficients $a_k(f)$ of $f$ are defined by
\begin{equation}\label{fouco}
a_k(f)=\int_Sf(y)R_k(y)\,d\mu(y).\end{equation} $s_n(f)$ denote
the partial sum of the generalized Fourier series of $f,$ i.e.
\begin{equation}\label{fou}
s_n(f,x)=\sum_{k=0}^n a_k(f)R_k(x)h(k).
\end{equation}
\begin{thm}
 Let $\{R_n\}_{n=0}^\infty$ be orthogonal polynomials satisfying (\ref{basic}-\ref{last1}).
 Then for any $f\in C(S)$ the partial sums $s_n(f,x)$ are convergent to $f$ uniformly on $S.$
\end{thm}
\begin{proof}
By orthogonality we have that
$s_n(R_m,x)=R_m(x) $ for $n\ge m.$ Therefore for any polynomial $p(x)$ there holds
$s_n(p,x)=p(x)$ for $n\ge \deg p.$ Since the polynomials are dense in $C(S)$ (as $S$ is a compact
subset of the real line) it suffices
to show that partial sums are uniformly bounded in $L^\infty$ norm, i.e. there exists a constant
$c$ such that
\begin{equation}\label{main}\|s_n(f,x)\|_{L^\infty}\le c\|f\|_{L^\infty}.
\end{equation}
The proof of this estimate will go roughly along the lines of \cite{O,OS1}, except that we have to overcome
technical
difficulties arising from the fact that orthogonality measure is not given explicitly.
By  (\ref{orth}) we get
$$s_n(f,x)=\int_S f(y)\sum_{k=0}^nR_k(x)R_k(y)h(k)\,d\mu(y)=\int_S f(y)\sum_{k=0}^np_k(x)p_k(y)\,d\mu(y).$$
Define the generalized Dirichlet kernel $K_n(x,y)$ by
\begin{equation}\label{dir}
K_n(x,y)=\sum_{k=0}^np_k(x)p_k(y).
\end{equation}
Then
\begin{multline*}
\|s_n(f,x)\|_{L^\infty}=\sup_{x\in S}\left |\int_S f(y)K_n(x,y)\,d\mu(y)\right |\\
\le \|f\|_{L^\infty}\sup_{x\in S}\int_S |K_n(x,y)|\,d\mu(y).
\end{multline*}
The proof will be finished if we show that
\begin{equation}\label{est}
\sup_{x\in S}\int_S |K_n(x,y)|\,d\mu(y)<+\infty.
\end{equation}
For this purpose we will use the conclusion of Theorem 1
which implies in particular that $S=\{0\}\cup\{\xi_k\}_{k=1}^\infty$ and
$\xi_n\sim q^n.$
Since $S\subset [0,\xi_1]$   we obtain
\begin{equation}\label{inte}
\int_S|K_n(x,y)|\,d\mu(y)=\int_{[0,\xi_n]}|K_n(x,y)|\,d\mu(y)+
\int_{(\xi_n,\xi_1]}|K_n(x,y)|\,d\mu(y).
\end{equation}
Combining (\ref{s}), (\ref{B}) and (\ref{hn}) yields
$$ \int_{[0,\xi_n]}|K_n(x,y)|\,d\mu(y)\le \mu([0,\xi_n])\sum_{k=0}^np_k^2(0)
\le c ,$$ for some constant independent of $n.$ It remains to
estimate uniformly the second integral of the right hand side of
(\ref{inte}) for $x\in S=\{0\}\cup \{\xi_k\}_{k=1}^\infty.$ We
split this integral into
$$K_n(x,x)\mu(x)+\int_{(\xi_n,\xi_1],y\neq x}|K_n(x,y)|d\mu(y).$$
The first term is less than 1, because
$$\mu(x)^{-1}=\sum_{k=0}^\infty p_k^2(x)\ge \sum_{k=0}^n p_k^2(x)=K_n(x,x).$$

By the  Christoffel-Darboux formula (\cite[1.17]{AH} we have
$$K_n(x,y)=\lambda_n{p_{n+1}(x)p_n(y)-p_n(x)p_{n+1}(y)\over x-y}.$$
Moreover since $\xi_{k+1}/\xi_k$ is bounded away from 1 there exists
a constant $d$ such that
$$|x-y|\ge dy,\quad x\neq y,\ x,y\in \{\xi_i\}_{i=1}^\infty.$$
Therefore by using $|p_k(x)|\le p_k(0)$ for $x\in S$ we obtain
\begin{multline*}
\int_{(\xi_n,\xi_1],y\neq x}|K_n(x,y)|d\mu(y)\le \\
{\lambda_np_{n+1}(0)\over d}\int_{(\xi_n,\xi_1]}{|p_n(y)|\over
y}\,d\mu(y)+ {\lambda_np_{n}(0)\over
d}\int_{(\xi_n,\xi_1]}{|p_{n+1}(y)|\over y}\,d\mu(y).
\end{multline*}
In view of $\lambda_n=\sqrt{\alpha_{n+1}\gamma_n}\approx q^n$ and $p_n(0)\sim s^{n/2}$
(see (\ref{s}) and (\ref{hn})) it suffices to show that
\begin{equation}\label{wazne}\int_{(\xi_n,\xi_1]}{|p_n(y)|\over y}\,d\mu(y)=O(q^{-n}s^{-n/2}).
\end{equation}
Fix  a nonnegative integer $l$ such that $q^{2l+2}<s^{-1}.$ Then we have
\begin{multline*}
\left(\int_{(\xi_n,\xi_1]}{|p_n(y)|\over y}\,d\mu(y)\right )^2=
\left(\int_{(\xi_n,\xi_1]}{y^l|p_n(y)|\over y^{l+1}}\,d\mu(y)\right )^2\\\le
\int_{S}y^{2l}p_n^2(y)\,d\mu(y)\ \int_{(\xi_n,\xi_1]}{1\over y^{2l+2}}\,d\mu(y).
\end{multline*}
Then we apply the recurrence relation (\ref{recur}) $2l$ times, and use orthonormality and the fact
that $\beta_n\approx q^n,$ $\lambda_n\approx q^n,$   to get
$$ \int_{S}y^{2l}p_n^2(y)\,d\mu(y)=O(q^{2nl}).$$
On the other hand by Theorem 1 we have $\xi_k\le C q^k$ and $\mu(\{\xi_k\})\le Cs^{-k}$ for some
constant $C.$ Thus
\begin{multline*}
\int_{(\xi_n,\xi_1]}{1\over y^{2l+2}}\,d\mu(y)=
\sum_{k=1}^{n-1}\xi_k^{-(2l+2)}\mu(\{\xi_k\})\\
\le C^2\sum_{k=1}^{n-1}q^{-k(2l+2)}s^{-k}=O(q^{-n(2l+2)}s^{-n}).
\end{multline*}
Therefore
$$\left(\int_{(\xi_n,\xi_1]}{|p_n(y)|\over y}\,d\mu(y)\right )^2=O(q^{-2n}s^{-n}),$$
as we required in (\ref{wazne}).
\end{proof}
{\bf Example.}
Fix $0<a<1$ and $0<q<1.$
Let $\alpha_n=a^2q^n$ and $\gamma_n=q^n.$
Then
$$
\beta_n=(1+a^2)q^n,\qquad \lambda_n=aq^{1/2}q^n.$$
It can be checked easily that the assumptions (\ref{compact})-(\ref{last1}) are satisfied
with $s=a^{-2},$ $\kappa=1,$ $N=1,$ i.e. there exists $c$ satisfying  (\ref{ad}) and (\ref{last1}), if
$${a\over 1+a^2}< q^{1/2}{1-q\over 1+q^2}.$$
Therefore for orthonormal polynomials associated with the recurrence
relation
$$xp_n=-\lambda_np_{n+1}+\beta_np_n-\lambda_{n-1}p_{n-1}$$
the conclusion of Theorem 2 holds. Moreover these polynomials
admit nonnegative product linearization.

\end{document}